\documentclass[11pt]{article}

\usepackage{amssymb,amsmath}
\usepackage{latexsym}
\usepackage{natbib,enumerate}

{\catcode `\@=11
\global\let \AddToReset=\@addtoreset}

\usepackage{subfigure}
\usepackage{graphicx}

\AddToReset{equation}{section}
\newtheorem{theorem}{Theorem}[section]
\newtheorem{lemma}{\bf Lemma}[section]

\newtheorem{proposition}{\bf Proposition}[section]

\newtheorem{@definition}{\sc Definition}[section]
\newtheorem{@remark}{\sc Remark}[section]
\newenvironment{remark}{\begin{@remark}\rm}{\end{@remark}}
\newtheorem{@example}{\sc Example}[section]
\newenvironment{example}{\begin{@example}\rm}{\end{@example}}
\def\mathsf{\bf}

\def\supp{\rm supp}

\def\E{\mathrm E}

\def\d{\mathrm d}
\def\e{\mathrm e}
\def\i{\mathrm i}

\def\text{\mbox}
\def\veps{\varepsilon}
\newcommand{\Var}{{\rm Var}}

\newcommand{\beq}{\begin{equation}}
\newcommand{\eeq}{\end{equation}}
\newcommand\beqn{\begin{displaymath}}  
\newcommand\eeqn{\end{displaymath}}
\newcommand\refeq[1]{{\rm (\ref{e:#1})}}

\reversemarginpar

\def\convd{\stackrel{\mbox{$\scriptstyle d$}}{\longrightarrow}}
\def\convP{\stackrel{\mbox{$\scriptstyle P$}}{\longrightarrow}}

\def\veps{\varepsilon}

\author{Dmitrij Celov${}^1$, Remigijus Leipus${}^{1,2}$
and Anne Philippe${}^3$ \bigskip \\
{\small ${}^1$Vilnius University, Lithuania}\\
{\small ${}^2$Institute of Mathematics and Informatics, Lithuania}\\
{\small ${}^3$Laboratoire de Math\'ematiques Jean Leray,
Universit\'e de
  Nantes, France}
}

\title{Asymptotic normality of the mixture density estimator in a disaggregation scheme\thanks{
The research was supported by bilateral France-Lithuanian
scientific project Gilibert and Lithuanian State Science and
Studies foundation (V-07058).}}
\date{\today}

\begin{document}
\maketitle{}

\vskip1cm

\begin{abstract}

The paper concerns the asymptotic distribution of the mixture
density estimator, proposed by
\cite{oppenheim:leipus:philippe:viano:2006}, in the
aggregation/disaggregation problem of random parameter AR(1)
process. We prove that, under mild conditions on the
(semiparametric) form of the mixture density, the estimator is
asymptotically normal. The proof is based on the limit theory for
the quadratic form in linear random variables developed by
\cite{Bhansali-etal2007}. The moving average representation of the
aggregated process is investigated. A small simulation study
illustrates the result.

 \vskip1cm


\vskip1cm

\end{abstract}
\newpage

\section{Introduction}

Aggregated time series data appears in different fields of studies
including applied problems in 
hydrology, sociology, statistics, economics. Considering
aggregation as a {\em time series} object, a number of important
questions arise. These comprise the properties of macro level data
obtained by small and large-scale aggregation in time, space or
both, assumptions of when and how the inverse (disaggregation)
problem can be solved, finally, how to apply theoretical  results
in practice.

Aggregated time series, in fact, can be viewed as a
transformation of the underlying time series by some (either
linear or non-linear) specific function defined at (in)finite
set of individual processes. In this paper we consider a linear
aggregation scheme, which is natural in applications. 
In practice it is found convenient to approximate individual data
by simple time series models, such as AR($1$), GARCH($1,1$) for
instance (see \cite{lewbel1994}, \cite{chong2006},
Zaffaroni~(\citeyear{zaffaroni2004}, \citeyear{zaffaroni2006})),
whereas
more complex individual data models do not provide an advantage in
accuracy and efficiency of estimates, and usually are very
difficult to study from the theoretical point of view.

Aggregation by appropriately averaging the micro level time series
models can give intriguing results. It was shown in Granger (1980)
that the large-scale aggregation of infinitely many short memory
AR(1) models with random coefficients can lead to a long memory
fractionally integrated process. It means that the properties of an
aggregate time series may in general differ from those of individual
data.

It is clear however that the weakest point of the aggregation is a
considerable loss of information about individual characteristics
of the underlying data. Roughly speaking, an aggregated time
series can not be as informative about the attributes of
individual data as the micro level processes are. On the other
hand, using some special aggregation schemes, which involve, for
instance, independent identically distributed ``elementary''
processes with known structure (such as AR($1$)), enables
to solve an inverse problem: to recover the properties of individual
series with the aggregated data at hand. This problem is called a
\textit{disaggregation problem}.

Different aspects of this problem were investigated in
Dacunha-Castelle, Oppenheim (2001),
\cite{oppenheim:leipus:philippe:viano:2006}, \cite{clp2007}. The
last two papers deal with the asymptotic statistical theory in the
disaggregation problem: they present the construction of the
mixture density estimate of the individual AR(1) models, the
consistency of an estimate, and some theoretical tools needed
here. Resuming the previous research, the major objective of the
present paper is to obtain the \textit{asymptotic normality}
property of the mixture density estimate, that enlarges the range
of applications, solving the accuracy of simulation studies,
statistical inference, forecasting and other problems.

Section~\ref{s:2} describes the disaggregation scheme, including
the construction of mixture density estimate proposed by
\cite{oppenheim:leipus:philippe:viano:2006}, and formulates the
main result of the paper. Important issues about the moving
average representation of the aggregated process are discussed in
Section~\ref{s:3}. The proof of the main theorem and auxiliary
results are given respectively in Section~\ref{s:proofmain} and
Section~\ref{s:app}. Some simulation results are presented in
Section~\ref{s:simul}.

\section{Preliminaries and the main result}\label{s:2}

Consider a sequence of independent processes $Y^{(j)}=\{Y^{(j)}_t,
t\in \bf Z\}$, $j\ge 1$ defined by the random coefficient AR(1)
dynamics \beq\label{e:ar1}
 Y_t^{(j)}=a^{(j)} Y^{(j)}_{t-1}+\varepsilon^{(j)}_t,
\eeq
 where $\veps^{(j)}_t, \ t\in {\bf Z}$, $ j=1,2,\dots$ are
 independent identically distributed (i.i.d.) random variables with $\E\veps^{(j)}_t=0$ and
 $0<\sigma_\veps^2=\E(\veps^{(j)}_t)^2<\infty$; $a, \, a^{(j)},\ j=1,2,\dots$ are i.i.d.\ random
 variables with $|a|\le 1$ and satisfying
\beq\label{e:a2}
 \E\Big[\frac{1}{1-a^2}\Big] <\infty.
\eeq
 It is assumed that the sequences $\{\veps^{(j)}_t, t\in {\bf Z}\}$, $j=1,2,\dots$ and
$\{a,a^{(j)}, j=1,2,\dots\}$ are independent.

Under these conditions, \refeq{ar1} admits a stationary solution
$Y^{(j)}$ and, according to \cite{oppenheim:viano:2004}, the
finite dimensional distributions of the process \beqn
 X^{(N)}_t=\frac{1}{\sqrt N}\; \sum_{j=1}^N Y_t^{(j)}, \ \ t\in {\bf Z},
\eeqn
 weakly converge as $N\to\infty$ to those of a zero mean stationary Gaussian
  process $X=\{X_t,t\in {\bf Z}\}$, called the {\em aggregated} process.
 Suppose that random coefficient $a$ admits a density $\varphi(x)$, absolutely continuous
 with respect to the Lebesgue measure, which by
  \refeq{a2} satisfies
\beq\label{e:mixd}
 \int_{-1}^1 \frac{\varphi(x)}{1-x^2}\ \d x<\infty.
\eeq
 Any density function satisfying \refeq{mixd} will be called a
 {\em mixture density}.

 Note that the covariance function and the spectral density of aggregated process
$X$ coincides with those of $Y^{(j)}$ and are given, respectively,
by
\beq\label{e:cov}
 \sigma(h):={\rm Cov}(X_h,X_0)=\sigma^2_\veps\int_{-1}^1 \frac{x^{|h|}}{1-x^2}\ \varphi(x)\d x
\eeq
 and
\beq\label{e:fl}
 f(\lambda) =\frac{\sigma^2_\veps}{2\pi} \int_{-1}^1 \frac{\varphi(x)}{|1-x\e^{\i\lambda}|^2} \
 \d x.
\eeq

 The {\em disaggregation problem} deals with finding the individual processes (if they exist) of
 form \refeq{ar1}, which produce the aggregated process
 $X$ with {\em given} spectral density $f(\lambda)$ (or covariance $\sigma(h)$). This is equivalent to finding  $\varphi(x)$ such that
 \refeq{fl} (or \refeq{cov}) and \refeq{mixd} hold.
 In this case, we say that {\em the mixture density $\varphi(x)$
 is associated with the spectral density $f(\lambda)$}.

 In order to estimate the mixture density
$\varphi(x)$ using aggregated observations $X_1,\dots,X_n$,
\cite{oppenheim:leipus:philippe:viano:2006} proposed the estimate
based on a decomposition of function
$\zeta(x)=\varphi(x)(1-x^2)^{-\alpha}$ in the orthonormal
$L^2(w^{(\alpha)})$--basis of Gegenbauer polynomials
$\{G_k^{(\alpha)}(x),k=0,1,\dots\}$, where
$w^{(\alpha)}(x)=(1-x^2)^\alpha$, $\alpha>-1$.  This decomposition
is valid (i.e.\ $\zeta$ belongs to $L^2(w^{(\alpha)})$) if
\beq\label{e:imp}
 \int_{-1}^1 \frac{\varphi^2(x)}{(1-x^2)^\alpha}\;\d x<\infty, \ \
 \alpha>-1.
\eeq

 Let $G_n^{(\alpha)}(x)=\sum_{j=0}^k g^{(\alpha)}_{n,j} x^j$.
The resulting estimate has the form
\beq\label{e:estimator}
 \hat\varphi_n(x)= \hat\sigma^{-2}_{n,\veps} (1-x^2)^\alpha \sum_{k=0}^{K_n}
 \hat\zeta_{n,k} G_k^{(\alpha)}(x),
\eeq
 where the $\hat\zeta_{n,k}$ are
estimates of the coefficients $\zeta_k$ in the $\alpha$-Gegenbauer
expansion of the function $\zeta(x)=\sum_{k=0}^\infty\zeta_k
G_k^{(\alpha)}(x)$ and are given by
 \beq\label{e:hz}
 \hat \zeta_{n,k} = \sum_{j=0}^k g_{k,j}^{(\alpha)}  (\hat\sigma_n(j)-
 \hat\sigma_n(j+2)),
\eeq
 $\hat\sigma^2_{n,\veps}=\hat
\sigma_n(0)-\hat\sigma_n(2)$ is the consistent estimator of
variance $\sigma^2_\veps$
 and $\hat\sigma_n(j)=n^{-1}\sum_{i=1}^{n-j} X_i X_{i+j}$ is the sample covariance of the
 aggregated process. Truncation level $K_n$ satisfies
\begin{equation}\label{eq:Kn}
    K_n = [\gamma \log n], \qquad 0<\gamma < (2\log
(1+\sqrt{2}))^{-1}.
\end{equation}

\cite{oppenheim:leipus:philippe:viano:2006} assumed the following
semiparametric form of the mixture density:
\beq\label{e:phi2}
 \varphi(x)=(1-x)^{1-2d_1}(1+x)^{1-2d_2} \psi(x), \ \ 0<d_1, d_2<1/2,
\eeq
 where $\psi(x)$ is continuous on $[-1,1]$ and does not vanish at
 $\pm 1$. Then, under conditions above and corresponding relations between
$\alpha$ and $d_1,d_2$, they showed the consistency of the
estimator $\hat\varphi_n(x)$ assuming that the variance of the
noise, $\sigma^2_\veps$, is {\em known} and equals 1. In more
realistic situation of unknown $\sigma^2_\veps$, it must be
consistently estimated. In order to understand intuitively the
construction of estimator $\hat\sigma^2_{n,\veps}$, it suffices to
note that, by \refeq{cov}, $\sigma^2_\veps=\sigma(0)-\sigma(2)$.
Also note that the estimator $\hat\varphi_n(x)$ in
\refeq{estimator} possesses property
$\int_{-1}^1\hat\varphi_n(x)\d x=1$, which can be easily verified
noting that $\int_{-1}^{1} (1-x^2)^\alpha G_k^{(\alpha)}(x) \; \d
x = (g^{(\alpha)}_{0,0})^{-1}$ if $k=0$, and $=0$ otherwise,
implying \beqn
  \int_{-1}^{1} (1-x^2)^\alpha \sum_{k=0}^{K_n}
 \hat\zeta_{n,k} G_k^{(\alpha)}(x)  {\rm d} x
 =\hat\zeta_{n,0}/g^{(\alpha)}_{0,0}=
  \hat \sigma_n(0) -\hat
\sigma_n(2) \eeqn
 by \refeq{hz}.

\smallskip


In this paper, we further study the properties of the proposed
mixture density estimator. In order to formulate the theorem about
the asymptotic normality of estimator $\hat \varphi_n(x)$, we will
assume that aggregated process $X_t$, $t\in {\bf Z}$ admits the
following linear representation.
\medskip

\noindent{\sc Assumption A}\ \ Assume that $X_t$, $t\in {\bf Z}$
is a linear sequence \beq\label{e:linear}
 X_t=\sum_{j=0}^\infty\psi_j Z_{t-j},
\eeq
 where the $Z_t$ are i.i.d.\ random variables with
 zero mean, finite fourth moment and the coefficients
 $\psi_j$ satisfy
\beq\label{e:as1}
 \psi_j\sim c j^{d-1},\  \ |\psi_j-\psi_{j+1}|=O(j^{d-2}), \ 0<d<1/2
\eeq
 with some constant $c\ne 0$.
\medskip

We also introduce the following condition on the mixture density
$\varphi(x)$.

\medskip

\noindent{\sc Assumption B}\ \ Assume that mixture density
$\varphi$ has a form \beq\label{e:phi1}
 \varphi(x)=(1-x)^{1-2d}\psi(x), \ \ 0<d<1/2,
\eeq
 where $\psi(x)$ is a
 nonnegative function with ${\supp} (\psi)\subset [-1,1]$, continuous at $x=1$, $\psi(1)\ne
 0$.
\medskip

Note that, omitting in \refeq{phi2} the factor responsible for the
seasonal part, we thus obtain the corresponding 'long memory'
spectral density with singularity at zero (but not necessary at
$\pm\pi$) and the corresponding behavior of the coefficients
 $\psi_j$ in linear representation \refeq{MArepr}.

\begin{theorem}\label{t:main}
Let $X_t, t\in{\bf Z}$ be the aggregated process satisfying
Assumption A and corresponding to the mixture density given by
Assumption B. Assume that \refeq{imp} holds, and $d$ and $\alpha$
satisfy the following condition
 \beq\label{e:alp}
  -1/2< \alpha <\frac52 -4d.
\eeq
 Let $K_n$ be given in (\ref{eq:Kn}) with $\gamma$ satisfying
\beq\label{e:cc}
 0<\gamma<(2\log(1+\sqrt
 2))^{-1}\Big(1-\max\Big\{\alpha+4d-\frac32,0\Big\}\Big).
\eeq
  Then for every fixed $x\in (-1,1)$, such that $\varphi(x)\ne 0$, it holds
\beq\label{e:asN}
 \frac{\hat\varphi_n(x)
 -\E\hat\varphi_n(x)}{\sqrt{\Var(\hat\varphi_n(x))}}\convd
  {\rm N}(0,1).
\eeq
\end{theorem}

Proof of the theorem is given in Section~\ref{s:proofmain}.

\begin{remark} Suppose that $\varphi(x)$ satisfies Assumption B.
Then assumption \refeq{imp} is equivalent to $\int_{-1}^1
\psi^2(x)(1+x)^{-\alpha}\d x<\infty$ and $\alpha<3-4d$. The last
inequality is implied by \refeq{alp}.

\end{remark}

\begin{example}\label{ex:1} Assume two mixture densities
\beq\label{e:FARIMAmd}
 \varphi(x;d) = C_1(d) x^{d-1} (1-x)^{1-2d} (1+x) {\bf 1}_{(0,1]}(x),
 \ 0<d<1/2,
\eeq
 where $C_1(d)=\frac{\Gamma(3-d)}{2\Gamma(d)\Gamma(2-2d)}$, and
\beq\label{e:bbb}
 \varphi_g(x;\kappa)= C_2(\kappa) |x|^\kappa {\bf 1}_{[-a_*,0]}(x), \ \ \kappa>0,
\eeq
 where $0<a_*<1$, $C_2(\kappa)=(\kappa+1)(a_*)^{-\kappa-1}$.

According to Dacunha-Castelle and Oppenheim (2001), the spectral
density corresponding to $\varphi(x;d)$ is FARIMA(0,$d$,0)
spectral density \beq\label{e:flambda1}
  f(\lambda;d)=\frac{1}{2\pi} \
  \Big(2\sin\frac{|\lambda|}{2}\Big)^{-2d}.
\eeq
 Also, since the support of $\varphi_g$ lies inside $(-1,1)$,
 the spectral density $g(\lambda;\kappa)$ corresponding to $\varphi_g(x;\kappa)$
 is analytic function (see Proposition 3.3 in \cite{clp2007}).

 Consider the spectral density given by
\beq\label{e:flambda}
  f(\lambda)=f(\lambda;d) g(\lambda;\kappa), \ \ \lambda\in [-\pi,\pi].
\eeq
 It can be shown that the mixture density
 $\varphi(x)$ associated with $f(\lambda)$ \refeq{flambda}
 is supported on $[-a_*,1]$, satisfies Assumption B with $\psi(x)$
 which is continuous function on $[-a_*,1]$ and at the neighborhood
 of zero satisfies $\psi(x)  =O(|x|^d)$. This implies the validity of condition
\refeq{imp} needed to obtain the corresponding $\alpha$-Gegenbauer
expansion. For the proof of this example and precise asymptotics
of $\psi(x)$ at zero see Appendix~A.

Finally, the aggregated process $X$, obtained using such mixture
density $\varphi(x)$, satisfies Assumption A by
Proposition~\ref{l:wold}, which shows that assumptions A and B are
satisfied under general 'aggregated' spectral density
$f(\lambda)=f(\lambda;d)g(\lambda)$, where $g(\lambda)$ is
analytic function on $[-\pi,\pi]$ and the associated mixture
density is supported on $[-a_*,0]$ with some $0<a_*<1$.
\end{example}

\begin{remark}
Note that the 'FARIMA mixture density' \refeq{FARIMAmd}, due to
factor $x^{d-1}$,  does not satisfy \refeq{imp} and a
"compensating" density such as $\varphi_g(x;\kappa)$ in
\refeq{bbb} is needed in order to obtain the needed integrability
in the neighborhood of zero. Obviously, for the same aim, other
mixture densities instead of $\varphi_g(x;\kappa)$ \refeq{bbb} can
be employed.
\end{remark}


\section{Moving average representation of the aggregated process}
\label{s:3}

 In order to obtain the asymptotic normality result in
 Theorem~\ref{t:main}, an important assumption is that the aggregated
 process admits a linear representation with coefficients decaying
at an appropriate rate (see \cite{Bhansali-etal2007}). The related
issues about the moving average representation of the aggregated
process are discussed in this section.

From the aggregating scheme follows that any aggregated process
admits an absolutely continuous
 spectral measure. If, in addition, its spectral density, say,
 $f(\lambda)$ satisfies
\beq\label{e:reg1}
 \int_{-\pi}^\pi\log f(\lambda)\d \lambda >-\infty,
\eeq
 then the function
\beqn
 h(z)=\exp\Big\{\frac{1}{4\pi} \int_{-\pi}^\pi\frac{\e^{\i\lambda}+z}{\e^{\i\lambda}-z}\; \log f(\lambda)\d \lambda\Big\},
 \ \ |z|<1,
\eeqn
 is an outer function from the Hardy space $H^2$, does not vanish for $|z|<1$ and
 $f(\lambda)=|h(\e^{\i\lambda})|^2$. Then, by the Wold decomposition theorem,
 corresponding process $X_t$ is purely nondeterministic
 and has the MA($\infty$) representation (see Anderson~(\citeyear{andersonbook71}, Ch.\ 7.6.3))
\beq \label{e:MArepr}
 X_t=\sum_{j=0}^\infty \psi_j Z_{t-j},
\eeq
 where the coefficients $\psi_j$ are defined from the expansion of normalized outer function
 $h(z)/h(0)$, $\sum_{j=0}^\infty\psi_j^2<\infty$, $\psi_0=1$, and
 $Z_t=X_t-\widehat X_t$, $t=0,1,\dots$ ($\widehat X_t$ is the optimal linear
 predictor of $X_t$) is the innovation process, which is zero mean, uncorrelated,
 with variance
\beq \label{e:sigma22}
 \sigma^2=2\pi\exp\Big\{\frac{1}{2\pi}\int_{-\pi}^\pi\log f(\lambda)\d\lambda
 \Big\}.
\eeq
 By construction, the aggregated process is Gaussian,
 implying that the innovations $Z_t$ are i.i.d.\ N$(0,\sigma^2)$ random variables.
\medskip

 Next we focus on the class of semiparametric mixture densities
 satisfying Assumption B. As it was mentioned earlier, this form
 is natural, in particular it covers the mixture densities $\varphi_1(x;d)$ and $\varphi(x)$ in Example~\ref{ex:1}.

\begin{proposition} Let the mixture density $\varphi(x)$ satisfies Assumption
B. Assume that either
\smallskip

\noindent (i) ${\supp}(\psi)= [-1,1]$ and
$\tilde\psi(x)\equiv\psi(x)(1+x)^{2\tilde d-1}$ is continuous at
$-1$ and $\tilde\psi(-1)\ne 0$ with some $0<\tilde d<1/2$,
\smallskip

\noindent or
\smallskip

\noindent (ii) ${\supp}(\psi)\subset [-a_*,1]$ with some
 $0<a_*<1$.
\smallskip

\noindent
 Then the
 aggregated process admits a moving average representation \refeq{MArepr},
 where the $Z_t$ are Gaussian i.i.d.\ random variables with zero mean and
 variance \refeq{sigma22}.
\end{proposition}


\noindent{\sc Proof.} ($i$) We have to verify that \refeq{reg1}
holds. Rewrite $\varphi(x)$ in the form
\beqn
 \varphi(x)=(1-x)^{1-2d}(1+x)^{1-2\tilde d}\tilde \psi(x).
\eeqn
Proposition 4.1 in \cite{clp2007} implies
\beqn
 f(\lambda)\sim C_1 |\lambda|^{-2d}, \quad |\lambda|\to 0,
\eeqn
 with $C_1>0$. Hence
$\log f(\lambda) \sim \log C_1 - C_2
 \log|1-\e^{\i\lambda}|$, $|\lambda|\to 0$, where $C_2>0$.
For any $\epsilon>0$ choose $0<\lambda_0\le \pi/3$, such that
\beqn
 - \frac{\log f(\lambda)-\log C_1}{C_2 \log|1-\e^{\i\lambda}|}-1\ge -\epsilon, \ \
 0<\lambda\le \lambda_0.
\eeqn
 Since $-\log|1-\e^{\i\lambda}|\ge 0$ for $0\le \lambda\le \pi/3$, we obtain
\beq\label{e:reg2}
 \int_0^{\lambda_0} \log f(\lambda)\d \lambda \ge\lambda_0\log C_1- C_2(1-\epsilon)\int_0^{\lambda_0} \log |1-\e^{\i\lambda}|
 \d \lambda>-\infty
\eeq
 using the well known fact that $\int_{0}^\pi
\log|1-\e^{\i\lambda}| \d \lambda =0$. Similarly,
\beq\label{e:reg3}
   \int_{\pi-\lambda_0}^\pi \log f(\lambda) \d \lambda > -\infty.
\eeq

When $\lambda\in [\lambda_0, \pi-\lambda_0]$, there exist
$0<L_1<L_2<\infty$ such that
$$
 L_1 \le \frac{1}{2\pi|1-x \e^{\i\lambda}|^2}\le L_2
$$
 uniformly in $x\in (-1,1)$. Thus, by \refeq{fl}, $L_1\le f(\lambda)\le L_2$
for any $\lambda\in [\lambda_0,\pi-\lambda_0]$, and therefore
\beq\label{e:reg4}
   \int_{\lambda_0}^{\pi-\lambda_0} \log f(\lambda) \d \lambda>-\infty.
\eeq
 \refeq{reg2}--\refeq{reg4} imply inequality \refeq{reg1}.
\smallskip

\noindent The proof in case ($ii$) is analogous to ($i$) and,
thus, is omitted. \hfill$\Box$

 \begin{lemma}\label{l:analytic}
If the spectral density $g(\lambda)$ of the aggregated process
$X_t,t\in {\bf Z}$ is
 analytic function on $[-\pi,\pi]$, then $X_t$ admits representation
\beqn
 X_t=\sum_{j=0}^{\infty} g_j Z_{t-j},
\eeqn
 where the $Z_t$ are i.i.d.\ Gaussian random variables with zero mean and
variance \beq\label{e:sigma2}
 \sigma_g^2=2\pi\exp\Big\{\frac{1}{2\pi}\int_{-\pi}^\pi\log
 g(\lambda)\d\lambda
 \Big\}
\eeq
 and the $g_j$ satisfy $|\sum_{j=0}^\infty g_j|<\infty$, $g_0=1$.
\end{lemma}

\noindent{\sc Proof.} From Proposition~3.3 in \cite{clp2007}
 it follows that there exists $0<a_*<1$ such that
\beq\label{e:ss}
 g(\lambda) = \frac{\sigma^2_\varepsilon}{2\pi}\int_{-a_*}^{a_*}
 \frac{\varphi_g(x)}{|1-x \e^{\i\lambda}|^2}\;  \d x.
\eeq
 For all $x\in[-a_*,a_*]$ and $\lambda\in[0,\pi]$ we have
$$
 \frac{1}{|1-x \e^{\i\lambda}|^2 } \geq C_3 > 0,
$$
where $C_3=C_3(a_*)$. This and \refeq{ss} imply $\int_{0}^\pi \log
g(\lambda) \d\lambda>-\infty$. Finally,
 $|\sum_{j=0}^\infty g_j|<\infty$ follows from representation
\beqn
 g(\lambda)= \frac{\sigma^2_g}{2\pi} \Big|\sum_{j=0}^\infty g_j \e^{\i j\lambda}\Big|^2
\eeqn
 and the assumption of analyticity of $g$.
\hfill$\Box$

\begin{proposition}\label{l:wold} Let $X_t,t\in {\bf Z}$ be an aggregated process with spectral
density \beq\label{e:flambda2}
 f(\lambda)=f(\lambda;d)g(\lambda),
\eeq
 where $f(\lambda;d)$ is FARIMA spectral density \refeq{flambda1}
and $g(\lambda)$ is analytic spectral density.
 Then:

 \smallskip

\noindent  (i) if mixture density $\varphi_g(x)$ associated with
$g(\lambda)$ satisfies ${\supp}(\varphi_g) \subset [-a_*,0]$ with
some $0<a_*<1$, then $\varphi(x)$, associated with $f(\lambda)$,
satisfies Assumption B.
\smallskip

 \noindent (ii) $X_t$ admits a linear representation
\refeq{MArepr},
 where the $Z_t$ are Gaussian i.i.d.\ random variables with zero mean and variance
\beqn
 \sigma^2=2\pi\exp\Big\{\frac{1}{2\pi}\int_{-\pi}^\pi\log f(\lambda)\d\lambda\Big\}=
 \exp\Big\{\frac{1}{2\pi}\int_{-\pi}^\pi\log g(\lambda)\d\lambda\Big\} =\frac{\sigma^2_g}{2\pi}
\eeqn and the coefficients $\psi_j$ satisfy 
\beq\label{e:phi_as}
 \psi_j \sim \frac{\sum_{k=0}^\infty g_k}{\Gamma(d)} \; j^{d-1}, \ \ |\psi_j-\psi_{j+1}|=O(j^{d-2}),
\eeq
 where $\psi_0=1$. (Here, the $g_k$ are given in
 Lemma~\ref{l:analytic}.)
\end{proposition}

\noindent{\sc Proof.} ($i$) By Corollary 3.1 in \cite{clp2007},
the mixture density associated with the "product" spectral density
\refeq{flambda2} exists and has a form
\beq \label{e:ppp}
 \varphi(x) =C_*^{-1}  \bigg(\varphi(x;d)
 \int_{-a_*}^{0}\frac{\varphi_g(y)\d y}{(1-xy)(1-y/x)}
+ \varphi_g(x) \int_{0}^{1}\frac{\varphi(y;d)\d y}{(1-xy) (1-y/x)}
\bigg), \eeq
with \beq\label{e:C*}
  C_*:=\int_0^1\bigg(\int_{-a_*}^0
 \frac{\varphi(x;d)\varphi_g(y)}{1-xy}\;\d y\bigg)\d x,
\eeq
 where $\varphi(x;d)$ is given in \refeq{FARIMAmd}
 and is associated with the spectral density $f(\lambda;d)$, and
$\varphi_g(x)$ is associated with the spectral density
$g(\lambda)$. Clearly, this implies that Assumption B is
satisfied.

\smallskip

($ii$) We have
 \beqn
 f(\lambda;d)=
  \frac{1}{2\pi}\bigg|\sum_{j=0}^\infty h_j \e^{\i j\lambda}\bigg|^2\
 \ {\rm with} \ \ h_j =\frac{\Gamma(j+d)}{\Gamma(j+1)\Gamma(d)}
\eeqn
 and, recall,
\beqn
 g(\lambda)=\frac{\sigma^2_g}{2\pi}\bigg|\sum_{j=0}^\infty g_j \e^{\i j\lambda}\bigg|^2,\
  \ \ \ \sum_{j=0}^\infty g^2_j <\infty
\eeqn
 since, by Lemma~\ref{l:analytic}, $\int_{-\pi}^\pi \log
 g(\lambda)\d\lambda>-\infty$.
 On the other hand, $\int_{-\pi}^\pi \log f(\lambda)\d\lambda>-\infty$ implies
\beqn
 f(\lambda)=\frac{1}{2\pi}\Big|\sum_{j=0}^\infty \tilde\psi_j \e^{\i j\lambda}\Big|^2, \
 \ \
 \sum_{j=0}^\infty\tilde \psi_j^2<\infty
\eeqn
 and, by uniqueness of the representation,
\beqn
 \tilde\psi_k= \frac{\sigma_g}{\sqrt{2\pi}}\sum_{j=0}^k h_{k-j} g_j.
\eeqn
 It easy to see that,
\beq\label{e:dct}
  \sum_{j=0}^k h_{k-j} g_j\sim h_k \sum_{j=0}^\infty g_j\sim C_4 k^{d-1},
\eeq
 where $C_4= \Gamma^{-1}(d) \sum_{j=0}^\infty g_j$. Indeed, taking into
 account that $h_k\sim \Gamma^{-1}(d) k^{d-1}$, we can write
\beqn
 \sum_{j=0}^k h_{k-j} g_j = \Gamma^{-1}(d)k^{d-1}\sum_{j=0}^\infty a_{k,j}
 g_j,
\eeqn
 where $a_{k,j}= h_{k-j} \Gamma(d) k^{1-d} {\bf 1}_{\{j\le k\}}\to
 1$ as $k\to\infty$ for each $j$. On the other hand, we have $|a_{k,j}|\le C
 (1+j)^{1-d}$ uniformly in $k$ and, since the $g_j$ decay exponentially fast,
 the sum $\sum_{j=0}^\infty (1+j)^{1-d} |g_j|$ converges and
 the dominated convergence theorem applies to obtain \refeq{dct}.

 Hence, we can write
\beqn
 f(\lambda)=\frac{\sigma^2_g}{(2\pi)^2}\Big|\sum_{j=0}^\infty \psi_j \e^{\i j\lambda}\Big|^2, \
 \ \psi_0=1,
\eeqn
 where $\psi_j=\tilde \psi_j\sqrt {2\pi}/\sigma_g\sim C_4 j^{d-1}$. Thus, representation \refeq{MArepr} and the
 first relation in \refeq{phi_as} follows.

Finally, in order to check the second relation in \refeq{phi_as},
it suffices to note that \beqn
 \psi_j-\psi_{j+1}=\sum_{i=0}^j (h_{j-i}-h_{j+1-i}) g_i -g_{j+1},
\eeqn
 where $h_j-h_{j+1}\sim C_5 j^{d-2}$ and the $g_j$ decay
 exponentially fast.
 \hfill$\Box$

\section{Proof of main result}\label{s:proofmain}

In order to prove Theorem~\ref{t:main}, we use the result of
\cite{Bhansali-etal2007}, who considered the following quadratic
form \beqn
 Q_{n,X}=\sum_{t,s=1}^n d_n(t-s) X_t X_s,
\eeqn
 where the  $X_t$ are linear sequences satisfying Assumption A and the function $d_n(k)$
 satisfies the following assumption.
\medskip

\noindent{\sc Assumption C}\ \  Suppose that \beqn
 d_n(k)=\int_{-\pi}^\pi \eta_n(\lambda)
 {\e}^{{\i}k\lambda}{\d}\lambda
\eeqn
 with some even real function
 $\eta_n(\lambda)$, such that, for some $-1<\beta<1$ and a sequence of
 constants $m_n\ge 0$, it holds
\beq \label{e:bgk1}
 |\eta_n(\lambda)|\le m_n |\lambda|^{-\beta}, \ \lambda\in
 [-\pi,\pi].
\eeq \smallskip

Denote by $E_n$ a matrix $(e_n(t-s))_{t,s=1,\dots,n}$, where
\begin{equation}\label{eq:defEn}
  e_n(t-s)=\int_{-\pi}^\pi \eta_n(\lambda) f(\lambda) \e^{\i\lambda(t-s)}\d \lambda
\end{equation}
and let $\|E_n\|^2=\sum_{t,s=1}^n e_n^2(t-s)$.

\begin{theorem} {\rm  [\cite{Bhansali-etal2007}]}\label{t:bgk}
 Suppose that
 assumptions A and C are satisfied. If
 $2d+\beta<1/2$ and
\beq\label{e:rnEn}
 r_n=o(\|E_n\|),
\eeq
 where
\beq\label{e:rn}
 r_n =
\begin{cases}
  m_n n^{\max(0,2d+\beta) } & $if$ \ \ \text{$2d + \beta \not = 0$, } \\
  m_n \log n  &  $if$ \ \ \text{$2d + \beta  = 0$, }
\end{cases}
\eeq
 then, as
 $n\to\infty$, it holds
\beqn
 \Var (Q_{n,X})\asymp \|E_n\|^2
\eeqn
 and
\beqn
 \frac{Q_{n,X}- \E Q_{n,X}}{\sqrt{\Var (Q_{n,X})}} \convd {\rm N}(0,1).
\eeqn (Here for $a_n,b_n\ge 0$, $a_n\asymp b_n$ means that $C_6
b_n\le a_n\le C_7 b_n$ for some $C_6,C_7>0$.)
\end{theorem}

\noindent{\sc Proof of Theorem~\ref{t:main}}. First of all, note
that \beqn
 \hat\sigma^2_{n,\veps}\convP \sigma^2_\veps,
\eeqn
 which easily follows using Theorem 3 in \cite{hosking1996}.
 Hence, to obtain convergence \refeq{asN}, we can replace
 the factor $\hat\sigma^2_{n,\veps}$ by $\sigma^2_{\veps}$ in the definition
 of $\hat\varphi_n(x)$. Without loss of generality assume that $\sigma^2_\veps=1$.

 Rewrite the estimate $\hat \varphi_n(x)$ in a form
\begin{align}
\hat\varphi_n(x) &= (1-x^2)^\alpha \sum_{k=0}^{K_n} \sum_{j=0}^k
 g_{k,j}^{(\alpha)}  (\hat\sigma_n(j)-
 \hat\sigma_n(j+2)) G_k^{(\alpha)}(x)\nonumber \\
&=
 (1-x^2)^\alpha \sum_{k=0}^{K_n} G_k^{(\alpha)}(x) \sum_{j=0}^k
g_{k,j}^{(\alpha)} \int_{-\pi}^\pi (\e^{\i\lambda j}-\e^{\i\lambda
(j+2)})
I_n(\lambda)\d\lambda\nonumber \\
&=  \int_{-\pi}^\pi \eta_n(\lambda;x)
I_n(\lambda)\d\lambda,\label{eq:2}
\end{align}
 where
\beq\label{e:etan}
 \eta_n(\lambda;x):= (1-x^2)^\alpha \sum_{k=0}^{K_n} G_k^{(\alpha)}(x)
 \sum_{j=0}^k g_{k,j}^{(\alpha)} (\e^{i\lambda j}-\e^{i\lambda (j+2)})
\eeq
 and $I_n(\lambda)=(2\pi n)^{-1} |\sum_{j=1}^{n} X_j \e^{\i j\lambda}|^2$, $\lambda\in [-\pi,\pi]$
 is the periodogram.

Now the proof follows from Assumption A and the results obtained
in Lemma~\ref{l:eta} and Lemma~\ref{lem:norm} below, which imply
that, under appropriate choice of $m_n$ and $\beta$, all the
assumptions in Theorem~\ref{t:bgk} are satisfied. In particular,
by Lemma~\ref{l:eta}, the following bound for the kernel
$\eta_n(\lambda;x)$ holds
\beq\label{e:aaa}
 |\eta_n(\lambda;x)|\leq m_n |\lambda|^{-\beta},
\eeq
 where
\beq\label{e:mb}
 m_n=C_8 n^{\gamma\log(1+\sqrt 2)}, \ \ \beta=
 \frac{\alpha}{2}-\frac{3}{4},
\eeq
 $C_8$ is a positive constant, depending on $x$ and $\alpha$.
 Clearly, \refeq{alp} implies that $-1<\beta\le \frac12 -2d<\frac12$ and $2d+\beta<\frac12$.

Consider the cases $2d+\beta\le 0$ or $0<2d + \beta<1/2$. In the
case $2d+\beta\le 0$, from \refeq{rn}, \refeq{mb} we obtain \beqn
 r_n=C_8\begin{cases}
 n^{\gamma\log(1+\sqrt 2)}& {\rm if}\ \  2d+\frac{\alpha}{2}-\frac{3}{4}<0,\\
 n^{\gamma\log(1+\sqrt 2)}\log n & {\rm if}\ \
 2d+\frac{\alpha}{2}-\frac{3}{4}=0.
 \end{cases}
\eeqn
 Hence, by Lemma~\ref{lem:norm}, $r_n\|E_n\|^{-1}\to 0$ because $\gamma\log(1+\sqrt 2)<1/2$.

Assume now $2d + \beta>0$. Then
\beqn
  r_n=C_8n^{\gamma\log(1+\sqrt 2)+2d+\frac{\alpha}{2}-\frac{3}{4}}
\eeqn and
 $r_n\|E_n\|^{-1}\to 0$ by \refeq{cc}. \hfill$\Box$
\medskip

The following lemma shows that the kernel $\eta_n(\lambda;x)$
given in \refeq{etan} satisfies inequality \refeq{aaa} with $m_n$
and $\beta$ given in \refeq{mb}.

\begin{lemma} \label{l:eta}
For quantity $\eta_n(\lambda;x)$ given in \refeq{etan} and for
 every fixed $x\in (-1,1)$, $0<|\lambda|<\pi$ it holds
\beqn
 |\eta_n(\lambda;x)|\leq C_9 n^{\gamma \log (1+\sqrt 2)}
 |\lambda|^{(3-2\alpha)/4}
 \begin{cases}
  (1-x^2)^{\alpha/2-1/4}  \ \ $if$\ \ \alpha>-1/2,\\
  (1-x^2)^{\alpha}\ \  $if$\ \ -1<\alpha<-1/2,
 \end{cases}
\eeqn
 where $C_9$ depends on $\alpha$, and  $\gamma$ is given in (\ref{eq:Kn}).
\end{lemma}

\begin{lemma} \label{lem:norm}
Assume that a mixture density $\varphi(x)$ satisfies condition
\refeq{imp} and let $K_n\to\infty$. Then for every $x\in (-1,1)$,
such that $\varphi(x)\ne 0$ it holds \beq\label{eq:1}
 \|E_n\|^2 \ge C_{10} n (1+o(1)),
\eeq
 where $C_{10}>0$ is positive constant depending on $\alpha$ and $x$.
\end{lemma}

Proof of these two lemmas are given in Appendix~B.

\section{A simulation study}\label{s:simul}

In order to gain further insight into the asymptotic normality
property of the mixture density estimator
 \refeq{estimator}, in this section we conduct a Monte-Carlo simulation
 study. Several examples are considered, which correspond to the mixture densities
 having different shapes  (here we do not pose a question which rigorous
 aggregating schemes lead to the latter).

The following two families of mixture densities \beqn
 \varphi(x)= w \varphi_1(x)+(1-w)\varphi_2(x), \ \ 0<w<1,
\eeqn are considered:
\begin{itemize}
    \item Beta-type mixture densities defined by
\begin{eqnarray*}
\varphi_1(x)&\propto& x^{p_1-1}(1-x)^{q_1-1} {\bf 1}_{[0,1]}(x), \
\
    p_1>0,\
    q_1>0,\\
   \varphi_2(x)&\propto& |x|^{p_2-1}(a_*+ x)^{q_2-1}  {\bf 1}_{[-a_*,0]}(x), \ \
    p_2>0,\
    q_2>0,\ 0<a_*<1;
\end{eqnarray*}
    \item mixed (Beta and Uniform)-type mixture densities defined by
    \begin{eqnarray*}
\varphi_1(x)&\propto& x^{p_3-1}(1-x)^{q_3-1} {\bf 1}_{[0,1]}(x), \
\    p_3>0,\
    q_3>0,\\
   \varphi_2(x)&=& 
   a_*^{-1}
   \ {\bf 1}_{[-a_*,0]}(x),\ \ 0<a_*<1.
\end{eqnarray*}
\end{itemize}





In order to construct the mixture density estimator, in the first
step, the parameters $K_n$ and $\alpha$ must be chosen.
Preliminary Monte-Carlo simulations showed that the estimator
$\hat \varphi_n(x)$ has the minimal mean integrated square error
(MISE) when the parameter $\alpha$ is chosen to be equal $1-2d$.
The justification of this interesting conjecture remains an open
problem. This rule also ensures that \refeq{alp} is satisfied. The
number of Gegenbauer polynomials $K_n$ is chosen according to
(\ref{eq:Kn}). Note that,
 by construction, the estimator $\hat\varphi_n(x)$ is not necessary
positive, though it integrates to one.


In Figure~1
, we present three graphs and corresponding box plots for the
mixture densities of the form above. Cases 1 and 2 correspond to
the Beta-type mixture densities, Case 3 corresponds to the mixed
(Beta and Uniform)-type mixture density. The parameter values are
presented in Table~\ref{tab:1}. The box plots are obtained by a
Monte-Carlo procedure based on $M =500$ independent replications
with sample size $n=1500$ and bandwidth $K_n=3$ (we aggregate
$N=5000$ i.i.d.\ AR(1) processes). Individual innovations
$\veps^{(j)}_t$ are i.i.d.\ ${\rm N}(0,1)$. Note that the mixture
density in Case 2 corresponds to Example~\ref{ex:1} with the
parameters $d=0.2$, $\kappa=0.1$ (in the sense of behavior at
zero).

\begin{figure}[h!]\label{fig:compar}
\center{\subfigure[Case 1]{\includegraphics[scale=.30]{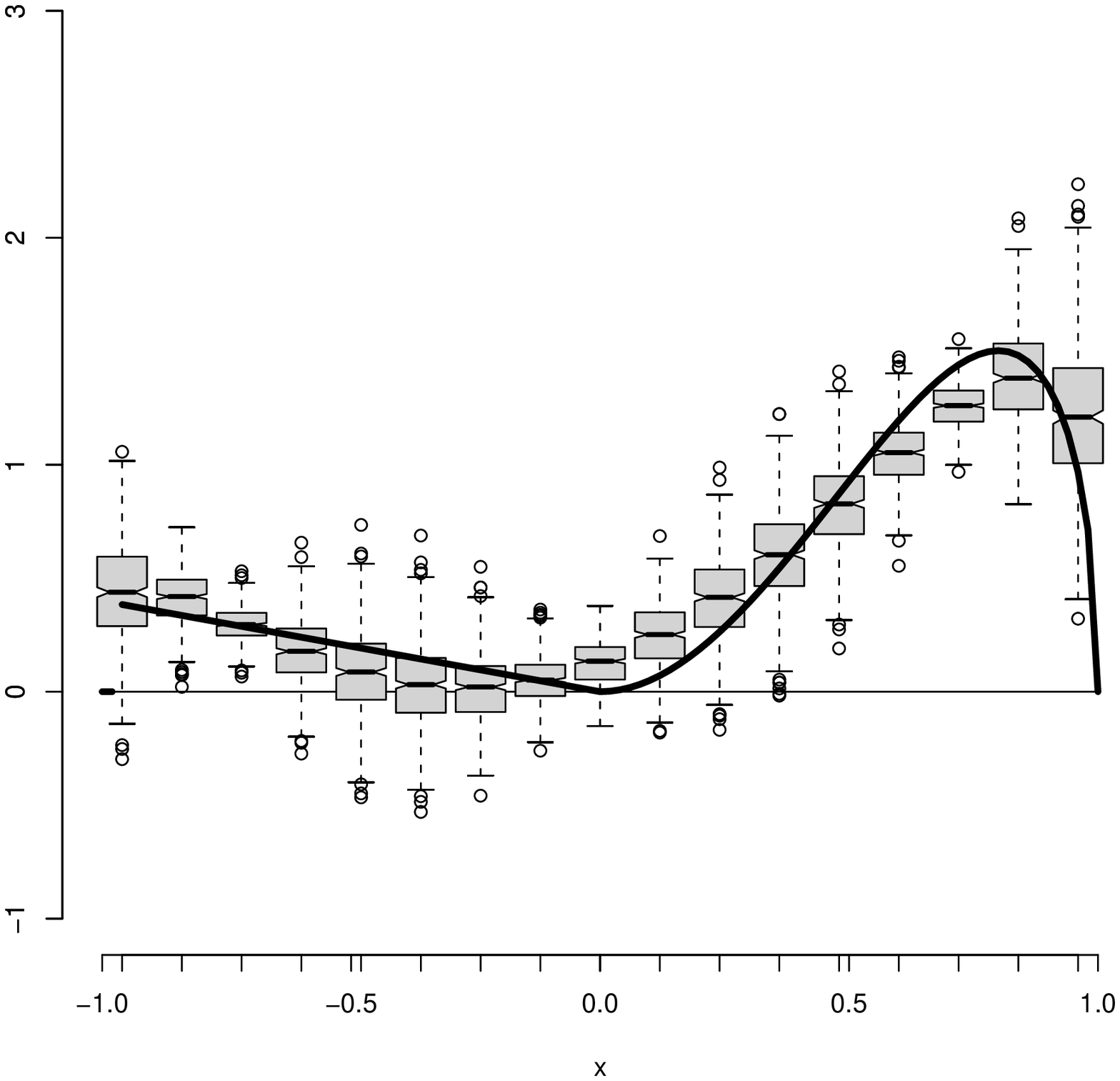}}
\subfigure[Case 2]{\includegraphics[scale=.30]{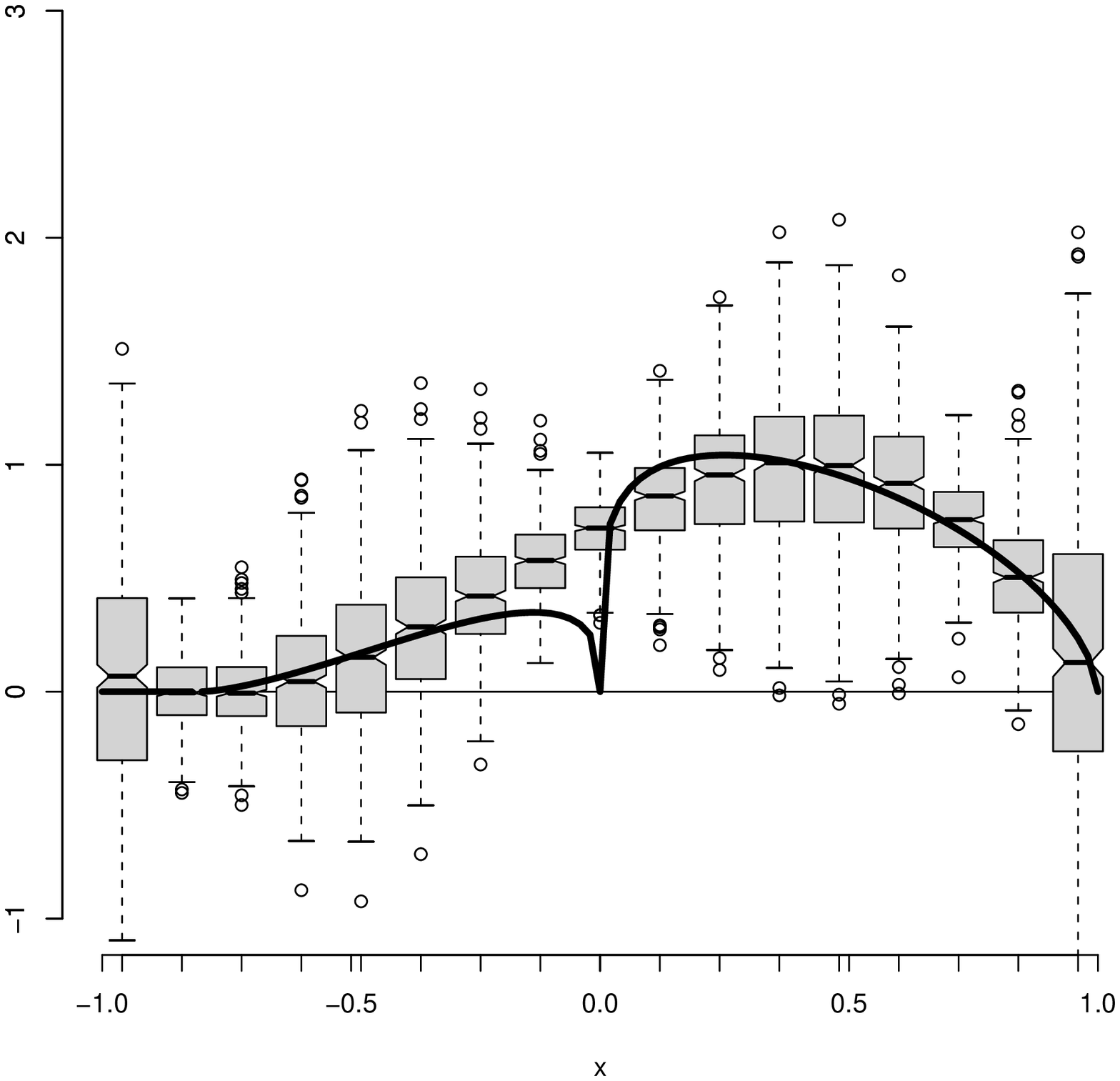}}}
\vspace{ -1cm}\\
\center{\subfigure[Case 3]{\includegraphics[scale=.30]{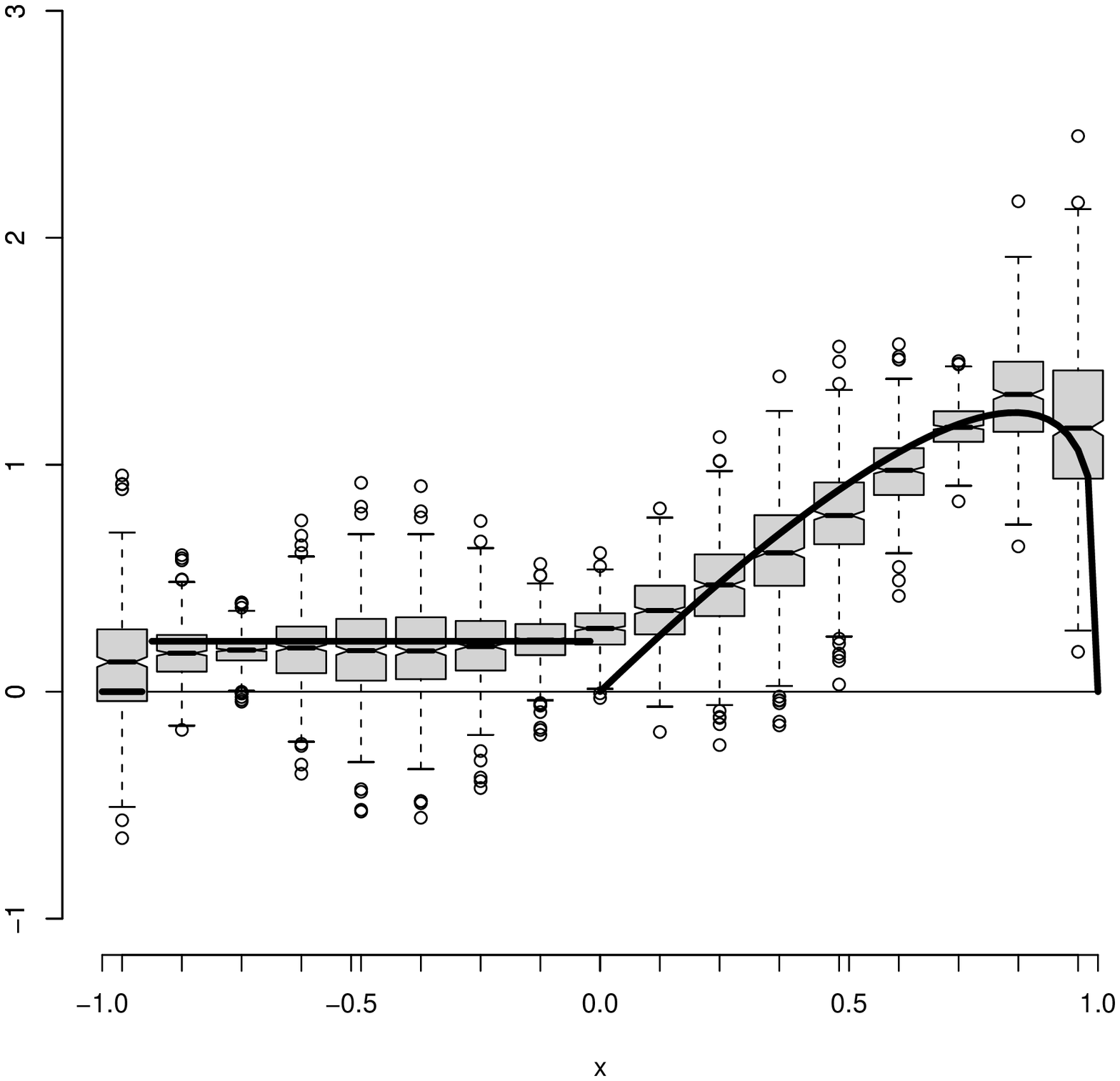}}}
\caption{True mixture densities (solid line) and the box plots of
the estimates. Number of replications $M=500$, sample size
$n=1500$.}
\end{figure}


\begin{table}[h]
    \centering
        \begin{tabular}[t]{c c c c c c c c}\hline
            { }    &  $w$& $a_*$ & $(p_1,q_1)$ & $(p_2,q_2)$ & $(p_3,q_3)$ & $d$ &  $\alpha$\\
            \hline
            Case 1&0.8&0.95&(3.0, 1.5)&(2.0, 1.0)& --& 0.25 & 0.5\\
            Case 2&0.8&0.80&(1.2, 1.6)&(1.3, 2.5)& --& 0.20 & 0.6 \\
            Case 3&0.8&0.90&-- & --& (2.0, 1.2)& 0.40 & 0.2\\
            \hline
        \end{tabular}
        \caption{Parameter values in cases 1--3.}\label{tab:1}
\end{table}

\begin{figure}[h!]\label{fig:2}
    \centering
        \includegraphics[width=12.5cm,height=9 cm
        ]{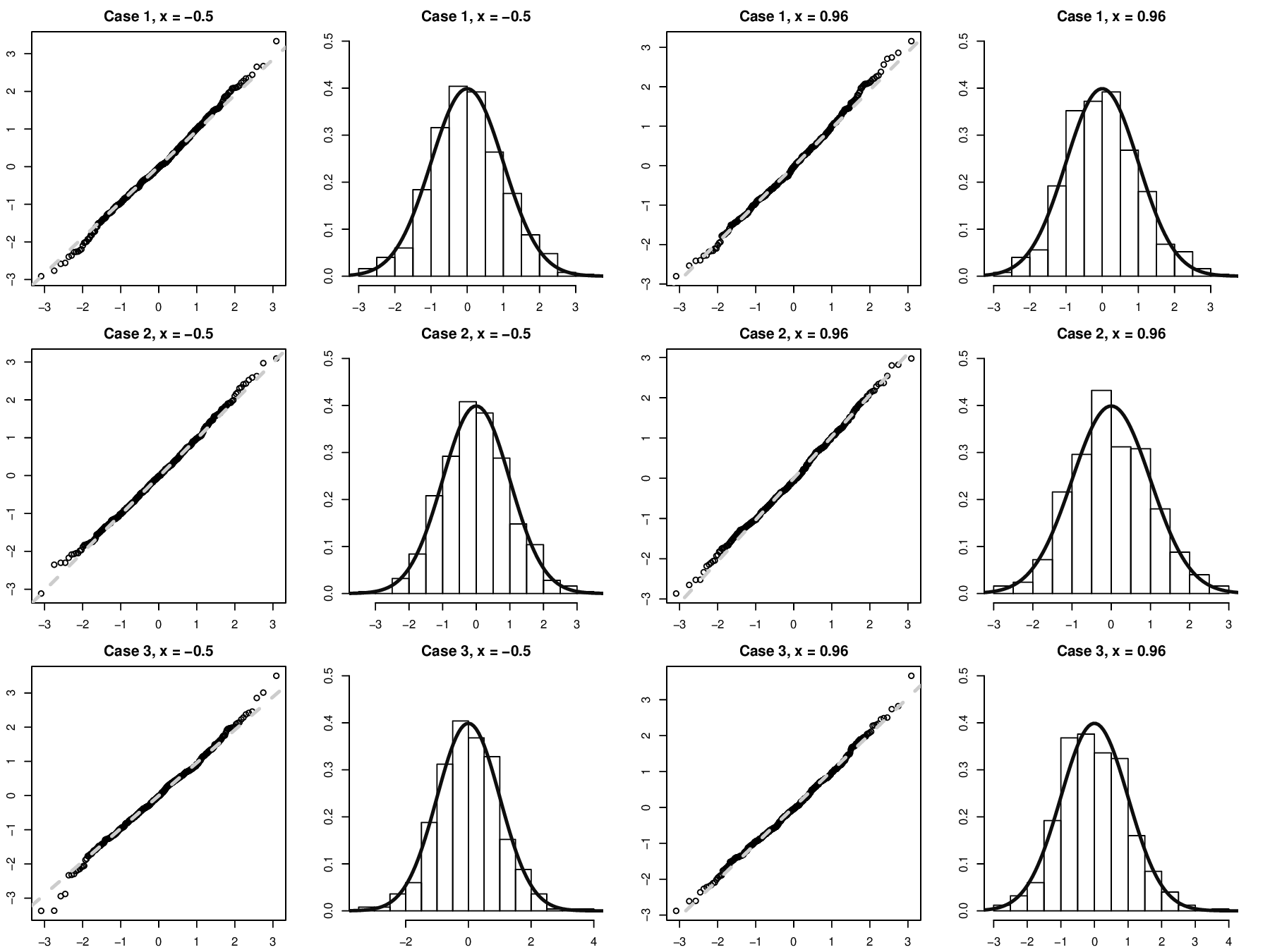}
    \caption{QQ plots and histograms of the estimates at points $x=-0.5$ and $x= 0.96$.
    Number of replications $M=500$, sample size
$n=1500$.}
\end{figure}

Box plots in Figure~1 
show that $\hat\varphi_n$ approximates the mixture density well
when $n$ is sufficiently large. However,  when the sample size is
relatively small it is difficult to estimate the mixture density
of the shape as in cases 2--3. This can be explained by the
construction of the estimator which assumes rather smooth form of
the mixture density around zero. On the other hand, it is clear
that the AR(1) parameter values which are close to zero does not
affect the long memory property. For our purposes, an important
fact is that the estimator correctly approximates the density at
the neighborhood of $x=1$. This enables us to estimate the unknown
(in real applications) parameter $d$ using a $\log$--$\log$
regression on periodogram at the neighborhood of this point (for
example Geweke and Porter-Hudak or Whittle-type estimators).

 Figure~2 supplements the earlier findings and shows that the distribution of estimator is approximately
 normal.\footnote{The Shapiro-Wilk test confirms
that
  in most cases normality hypothesis is consistent with the data.}
QQ-plots and histograms are given for fixed values $x=-0.5$ and
$x=0.96$ correspondingly. We use the same number of replications
$M=500$ and sample size $n=1500$.

\begin{figure}[h!]
\subfigure[Case 1, $x =
-0.5$]{\includegraphics[scale=.3]{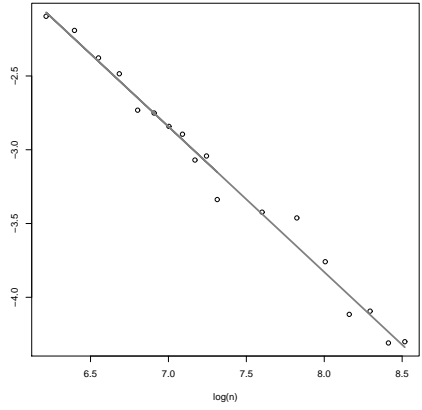}} \subfigure[Case 1, $x
= 0.96$]{\includegraphics[scale=.3]{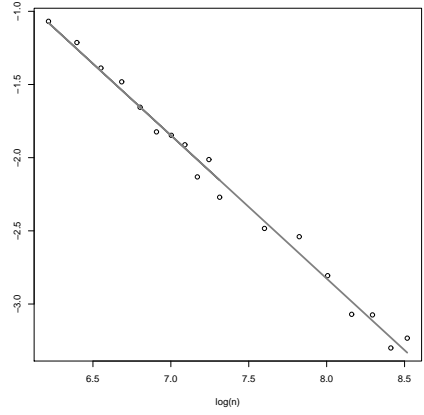}}
\caption{$\log$-$\log$ scale regression of the variance of
$\hat\varphi_n(x)$ as a function of $n$. The variance is estimated
using $M=500$ independent replications.}\label{fig:bpl}
\end{figure}

 The last Monte-Carlo experiment aims to show that the decay rate of $\Var
 (\hat{\varphi}_n(x))$ is $n^{-\gamma}$ with $\gamma=1$.
This ensures that the variance is decreasing fast enough.
To do this, we calculate
the $\log$--$\log$ regression of variance on the length
 of time series $n \in \{500,600,\dots,1400,1500,2000,\dots,5000\}$.
Figure~\ref{fig:bpl} demonstrates the corresponding parameter
estimates at different points and shows that $\hat{\gamma} \approx
1$.






\section{Appendix A. Proof of Example~\ref{ex:1}}

 By Corollary 3.1 in \cite{clp2007}, the mixture density
 $\varphi(x)$, $x\in [-a_*,1]$ associated with $f(\lambda)$ \refeq{flambda}
 is given by equality \refeq{ppp}, where $\varphi_g(x)\equiv
 \varphi_g(x;\kappa)$.
 Clearly, in this case, \refeq{ppp} can be rewritten in form \refeq{phi1} with
\beq\label{e:psipsi}
 \psi(x)= \tilde C (\psi_1(x)+\psi_2(x)),
\eeq
 where $\tilde C=C_1(d)C_2(\kappa)  C^{-1}_*$ is positive constant,
\begin{eqnarray}
 \psi_1(x)&:=&x^{d-1}(1+x){\bf 1}_{(0,1]}(x) \int_{-a_*}^0
 \frac{|y|^\kappa}{(1-xy)(1-y/x)}\;\d y,\label{e:psi1}\\
 \psi_2(x)&:=&|x|^\kappa(1-x)^{2d-1}{\bf 1}_{[-a_*,0]}(x) \int_0^1
 \frac{y^{d-1}(1-y)^{1-2d}(1+y)}{(1-xy)(1-y/x)}\;\d
 y.\label{e:psi2}
\end{eqnarray}
 Denote by $F(a,b;c;x)$ a hypergeometric function
$$
 F(a,b;c;x)=\frac{\Gamma(c)}{\Gamma(b)\Gamma(c-b)}\int_0^1
 t^{b-1}(1-t)^{c-b-1}(1-tx)^{-a} \d t,
$$
 with $c>b>0$ if $x<1$ and, in addition, $c-a-b>0$ if $x=1$. Then
 the corresponding integrals in $\psi_1(x)$ and $\psi_2(x)$ can be
 rewritten as
\begin{eqnarray*} &&\hspace{-1cm}\int_{-a_*}^0
 \frac{|y|^\kappa}{(1-xy)(1-y/x)}\;\d y,\\
&=&\frac{a_*^{\kappa+1}}{\kappa+1}
\frac{x(F(1,\kappa+1;\kappa+2;-a_*x)-
 F(1,\kappa+1;\kappa+2;-a_*/x))}{1-x^2}\\
 &\sim& \frac{a_*^{\kappa+1}}{\kappa+1} \ x, \ \ {\rm as}\ \ x\to 0+,
\end{eqnarray*}
and
\begin{eqnarray*}
 &&\hspace{-1cm}\int_0^1
 \frac{y^{d-1}(1-y)^{1-2d}(1+y)}{(1-xy)(1-y/x)}\;\d y\\
&=&\frac{\Gamma(d)\Gamma(2-2d)}{\Gamma(2-2d)}\frac{F(1,d;2-d;1/x)-x
 F(1,d;2-d;x)}{1-x}\\
&\sim& \Gamma(d)\Gamma(1-d) |x|^d, \ \ {\rm as}\ \ x\to 0-,
\end{eqnarray*}
where the last asymptotics follow from the well known properties
of the hypergeometric functions (see \cite{abramovitzstegun65}).

Thus, from \refeq{psi1}--\refeq{psi2} we obtain that
\begin{eqnarray}\label{e:qq}
 \psi_1(x) &\sim&\frac{a_*^{\kappa+1}}{\kappa+1} \ x^d, \ \ {\rm as}\ \ x\to 0+,\label{e:qq1}\\
 \psi_2(x)&\sim& \Gamma(d)\Gamma(1-d) |x|^{\kappa+d}, \ \ {\rm as}\ \ x\to 0-.\label{e:qq2}
\end{eqnarray}
\refeq{psipsi} and relations \refeq{qq1}--\refeq{qq2} complete the
proof.\hfill$\Box$

\section{Appendix B. Proofs of lemmas~\ref{l:eta}--\ref{lem:norm}}
 \label{s:app}

\noindent {\sc Proof of Lemma \ref{l:eta}.} By \refeq{etan},
\begin{eqnarray*}
  (1-x^2)^{-\alpha}\eta_n(\lambda;x)&=&\sum_{k=0}^{K_n} G_k^{(\alpha)}(x)
\sum_{j=0}^k g_{k,j}^{(\alpha)} (\e^{\i\lambda j}-\e^{\i\lambda (j+2)})\\
&=&(1-\e^{2\i\lambda})\sum_{k=0}^{K_n} G_k^{(\alpha)}(x)
\sum_{j=0}^k g_{k,j}^{(\alpha)} \e^{\i\lambda j}\\*
&=&(1-\e^{2\i\lambda})\sum_{k=0}^{K_n} G_k^{(\alpha)}(x)
 G_k^{(\alpha)}(\e^{\i\lambda}).
\end{eqnarray*}
 This and Lemma~\ref{lem:coeff} below implies
\begin{eqnarray}\label{e:sss}
  (1-x^2)^{-\alpha}|\eta_n(\lambda;x)|& \le & C_{11} |\lambda|^{-(2\alpha-3)/4}  \sum_{k=0}^{K_n} |G_k^{(\alpha)}(x)| (1+\sqrt 2)^k.
\end{eqnarray}
 Now, using the fact that for all $-1<x<1$
\beqn
 |G_k^{(\alpha)}(x)|\le \begin{cases}
  C_{12}(1-x^2)^{-\frac{\alpha}{2}-\frac{1}{4}} & {\rm if}\ \ \alpha>-1/2\\
                                C_{12} & {\rm if}\ \ \alpha<-1/2, \ \alpha\neq -3/2,-5/2,\dots
\end{cases}
\eeqn
(see inequality (7.33.6) in \cite{Szego} and (3.9) in
\cite{oppenheim:leipus:philippe:viano:2006})  and (\ref{eq:Kn}),
we get from \refeq{sss}
\begin{eqnarray*}
  (1-x^2)^{-\alpha}|\eta_n(\lambda;x)|&
  \le &  C_{13} |\lambda|^{-(2\alpha-3)/4}(1+\sqrt 2)^{K_n}\\
 & =&  C_{13}
 |\lambda|^{-(2\alpha-3)/4} \e^{K_n \log (1+\sqrt 2)} \\
  & \le &  C_9
 |\lambda|^{-(2\alpha-3)/4} n^{\gamma \log (1+\sqrt 2)}.
\end{eqnarray*}
\hfill$\Box$

\begin{lemma}\label{lem:coeff} For all $k\ge 0$, $\alpha>-1$, ($\alpha\ne -1/2$) and
$0<|\lambda|<\pi$ it holds
$$
  |(1-\e^{2\i\lambda})G_k^{(\alpha)}(\e^{\i\lambda})|\le C_{11} (1+\sqrt 2)^k
 |\lambda|^{-(2\alpha-3)/4},
$$
 where  constant $C_{11}$ depends on $\alpha$.
\end{lemma}

\noindent {\sc Proof.} Theorem 8.21.10 of \cite{Szego} implies
that for the usual (nonnormalized) Gegenbauer polynomials with
$\alpha>-1$, $\alpha\ne -1/2$ it holds
\begin{equation}
  \label{e:ck}
 C_k^{(\alpha+1/2)}(\e^{\i\lambda})=  \frac{\Gamma(k+\alpha+\frac12)}
{\Gamma(k+1)\Gamma(\alpha+\frac12)}\; z^k (1-z^{-2})^{-\alpha-1/2}
+ O(k^{\alpha-3/2} |z|^k),
\end{equation}
 where the complex numbers $w=\e^{\i\lambda}$ and $z$ are connected by
the elementary conformal mapping
\beq\label{e:conf}
 w=\frac{1}{2}(z+z^{-1}), \quad z= w+(w^2-1)^{1/2},
\eeq
 and $z$ satisfies $|z|>1$ (thus, $\lambda\ne 0,\pm\pi$).

Recall that the normalized Gegenbauer polynomials
$G_k^{(\alpha)}(z)$ are linked to $C_k^{(\alpha+1/2)}(z)$
 by equality
$$
 G_k^{(\alpha)}(z)=\gamma_k^{-1/2} C_k^{(\alpha+1/2)}(z), \ \textrm{where}\
 \gamma_k=\frac{\pi}{2^{2\alpha}} \frac{\Gamma(k+2\alpha+1)}
 {(k+\alpha+\frac{1}{2})\Gamma^2(\alpha+\frac{1}{2})\Gamma(k+1)}.
$$
Therefore, in terms of the normalized Gegenbauer polynomials, \refeq{ck}
 reads as follows
\begin{equation}\label{e:ck1}
 G_k^{(\alpha)}(\e^{\i\lambda})=  \frac{{\rm sgn}(\alpha+1/2) 2^\alpha}{\pi^{1/2}} \ b_k
 z^k (1-z^{-2})^{-\alpha-1/2} +O(k^{-1} |z|^k),
\end{equation}
where
$$
 b_k= \frac{(k+\alpha+1/2)^{1/2}\Gamma(k+\alpha+1/2)}
{\Gamma^{1/2}(k+1)\Gamma^{1/2}(k+2\alpha+1)}\to 1 \ \ {\rm as}\ \
k\to\infty.
$$

From \refeq{conf} we obtain for $w=\e^{\i\lambda}$
$$
w^2-1=\frac{1}{4} z^2(1-z^{-2})^2,
$$
which together with \refeq{ck1} yields
\begin{eqnarray*}
(1-\e^{2\i\lambda}) G_k^{(\alpha)} (\e^{\i\lambda})= -\frac{{\rm
sgn}(\alpha+1/2) 2^\alpha}{4\pi^{1/2}} \ b_k
 z^{k+2} (1-z^{-2})^{-\alpha+3/2} +O(k^{-1} |z|^k).
\end{eqnarray*}
Since $|z|>1$ and $z^2-1=2(\e^{2\i\lambda}-1)+2\e^{3\i\lambda/2}
 (\e^{\i\lambda}-\e^{-\i\lambda})^{1/2}$,
 we have
\begin{eqnarray}
 |1-z^{-2}|&\le& |z^2-1| \nonumber\\
 &\le& 2|\e^{2\i\lambda}-1|+2|\e^{\i \lambda}- \e^{-\i \lambda}|^{1/2}
\nonumber\\
 &=& 4|\sin\lambda|+2\sqrt 2 |\sin\lambda|^{1/2}\nonumber\\
 &\le& (4+2\sqrt 2) |\lambda|^{1/2}.\label{e:ck2}
\end{eqnarray}
So that, by \refeq{ck1}--\refeq{ck2},
\beq\label{e:O}
 |(1-\e^{2\i\lambda}) G_k^{(\alpha)} (\e^{\i\lambda})|\le
 C_{14} b_k |z|^k |\lambda|^{-(2\alpha-3)/4},
\eeq
 where $C_{14}=C_{14}(\alpha)$.

Finally, the straightforward verification shows that
\beqn
 \sup_{\lambda\in [-\pi,\pi]} |\e^{\i\lambda}+(\e^{2\i\lambda}-1)^{1/2}|=1+\sqrt 2.
\eeqn
 This completes the proof of lemma.\hfill $\Box$

\medskip

\noindent {\sc Proof of Lemma~\ref{lem:norm}}. Using
({\ref{eq:defEn}), \refeq{etan} rewrite the coefficients of $E_n$
\beqn
 e_n(t-s) = (1-x^2)^\alpha\sum_{k=0}^{K_n} G_k^{(\alpha)}(x) \sum_{j=0}^k g_{k,j}^{(\alpha)}
 \int_{-\pi}^\pi f(\lambda) (\e^{\i\lambda(t-s+j)}-\e^{\i\lambda(t-s+j+2)}) \d \lambda.
\eeqn Using the expression of the covariance function of an
aggregated process, we have for
 $t-s+j\ge 0$
\begin{eqnarray*}
\int_{-\pi}^\pi f(\lambda) (\e^{\i\lambda(t-s+j)}-\e^{\i\lambda(t-s+j+2)})\d \lambda
&=& \sigma(t-s+j)-\sigma(t-s+j+2) \\
&=& \sigma_\veps^2\int_{-1}^1 y^{t-s+j} \varphi(y) \d y.
\end{eqnarray*}
Thus, assuming $\sigma_\veps^2=1$, for $t-s\ge 0$ we have
\begin{eqnarray*}
e_n(t-s) &=& (1-x^2)^\alpha \sum_{k=0}^{K_n} G_k^{(\alpha)}(x)
\sum_{j=0}^k g_{k,j}^{(\alpha)} \int_{-1}^1 y^{t-s+j} \varphi(y) \d y  \\
&=& (1-x^2)^\alpha \sum_{k=0}^{K_n} G_k^{(\alpha)}(x)
\int_{-1}^1 y^{t-s}\varphi(y) \sum_{j=0}^k g_{k,j}^{(\alpha)}  y^{j}  \d y   \\
&=&  (1-x^2)^\alpha \sum_{k=0}^{K_n} G_k^{(\alpha)}(x)
\int_{-1}^1 y^{t-s}\varphi(y) G_{k}^{(\alpha)}(y) \d y.   \\
\end{eqnarray*}
 Integral $\int_{-1}^1 y^m\varphi(y) G_{k}^{(\alpha)}(y) \d y$ ($m$ is a nonnegative integer),
 appearing
 in the last expression is nothing else but the $k$th coefficient,
 $\psi_{m,k}$, in the $\alpha$-Gegenbauer expansion of the function
\beq\label{e:psim}
 \psi_{m}(x) =\frac{x^{m}\varphi(x)}{(1-x^2)^\alpha},
\eeq
 which obviously satisfies $\psi_{m}\in L^2(w^{(\alpha)})$. Therefore,
\begin{eqnarray*}
 e_n(t-s) &=& (1-x^2)^\alpha \sum_{k=0}^{K_n} G_k^{(\alpha)}(x) \psi_{|t-s|,k} \\
 &=& (1-x^2)^\alpha \bigg(\psi_{|t-s|}(x) - \sum_{k=K_n+1}^\infty G_k^{(\alpha)}(x) \psi_{|t-s|,k}\bigg)
\end{eqnarray*}
 and, denoting $R_n(m):=\sum_{k=K_n+1}^\infty G_k^{(\alpha)}(x) \psi_{|m|,k}$, $|m|<n$, we have
\begin{eqnarray*}
(1-x^2)^{-2\alpha}\|E_n\|^2 
 &=& \sum_{|m|<n} (n-|m|) \bigg(\psi_{|m|}(x)-\sum_{k=K_n+1}^\infty
 G_k^{(\alpha)}(x) \psi_{m,k}\bigg)^2\\
 &=& \sum_{|m|<n} (n-|m|) \psi^2_{|m|}(x) -2 \sum_{|m|< n} (n-|m|)
 \psi_{|m|}(x) R_n(m) \\
 &&+ \sum_{|m|< n} (n-|m|) R^2_n(m)=: A_{1,n} -2  A_{2,n} + A_{3,n}.
\end{eqnarray*}

Now, we prove that, as $n\to\infty$,
\beq\label{e:A1n}
 A_{1,n} \sim C_{15} n,
\eeq
 where $C_{15}=C_{15}(x)>0$ is some positive constant, and
\beq\label{e:A2n}
 A_{2,n} = o(n).
\eeq
  Since the last term $A_{3,n}$ is nonnegative by
construction, this will prove (\ref{eq:1}).

At points $x$ where $\varphi(x)\ne 0$ we have
\begin{eqnarray*}
 A_{1,n} 
 &=& \frac{\varphi^2(x)}{(1-x^2)^{2\alpha}} \sum_{|m|< n} (n-|m|) x^{2|m|}\\
 &\sim&  n\; \frac{\varphi^2(x)(1+x^2)}{(1-x^2)^{2\alpha+1}}, \ \
 {\rm as}\ \ n\to\infty,
\end{eqnarray*}
 which gives \refeq{A1n}. Consider term $A_{2,n}$. By \refeq{psim},
\begin{eqnarray*}
 A_{2,n}&=& \sum_{|m|<n} (n-|m|) \psi_{|m|}(x) \sum_{k=K_n+1}^\infty
 G_k^{(\alpha)}(x) \psi_{|m|,k} \\
 &= & \frac{\varphi(x)}{(1-x^2)^\alpha} \sum_{k=K_n+1}^\infty
 G_j^{(\alpha)}(x) \int_{-1}^1 \varphi(y) G_k^{(\alpha)}(y) \sum_{|m|< n} (n-|m|) (xy)^{|m|}
\d y \\
&=& \frac{\varphi(x)}{(1-x^2)^\alpha}\; (B_{1,n}-B_{2,n}-B_{3,n}),
\end{eqnarray*}
 where
\begin{eqnarray*}
 B_{1,n}&:=& n\sum_{k=K_n+1}^\infty
 G_j^{(\alpha)}(x) \int_{-1}^1 \varphi(y) G_k^{(\alpha)}(y) \sum_{m=-\infty}^\infty
 (xy)^{|m|}\d y,\\
 &=& n\sum_{k=K_n+1}^\infty
 G_j^{(\alpha)}(x) \int_{-1}^1 \varphi(y) G_k^{(\alpha)}(y) \frac{1+xy}{1-xy}\;\d y,\\
B_{2,n}&:=& \sum_{k=K_n+1}^\infty
 G_j^{(\alpha)}(x) \int_{-1}^1 \varphi(y) G_k^{(\alpha)}(y) \sum_{|m|< n}|m|
 (xy)^{|m|}\d y,\\
B_{3,n}&:=& n\sum_{k=K_n+1}^\infty
 G_j^{(\alpha)}(x) \int_{-1}^1 \varphi(y) G_k^{(\alpha)}(y) \sum_{|m|\ge n}
 (xy)^{|m|}\d y.
\end{eqnarray*}
 Since, by \refeq{imp},
\beqn
 \tilde \varphi_x(y) \equiv \frac{\varphi(y)}{(1-y^2)^\alpha}\frac{1+xy}{1-xy}
\eeqn
 satisfies $\tilde\varphi_x\in L^2(w^{(\alpha)})$ and $K_n\to\infty$, the sum
 $\sum_{k=K_n+1}^\infty$ in $B_{1,n}$ vanishes (as the tail of the convergent
 series). So that, $B_{1,n}=o(n)$ and, similarly, $B_{3,n}=o(n)$.

Finally,
\begin{eqnarray*}
 B_{2,n}&\sim& \sum_{k=K_n+1}^\infty
 G_j^{(\alpha)}(x) \int_{-1}^1 \varphi(y) G_k^{(\alpha)}(y)\;
 \frac{2xy}{(1-xy)^2}\;\d y=o(1)
\end{eqnarray*}
 using the similar argument as in the case of term $B_{1,n}$. This
 completes the proof of \refeq{A2n} and of the lemma. \hfill $\Box$

\bibliographystyle{apalike}

\end{document}